\documentclass[english]{paper}

\usepackage[foot]{amsaddr}

\usepackage{mama}
\usepackage{quiver}

\title{Seeing double through dependent optics}
\author{Matteo Capucci}
\email{matteo.capucci@strath.ac.uk}
\address{Office 1310, Livingstone Tower, 26 Richmond St, Glasgow (UK)}

\usepackage{amsbsy}
\usepackage{todonotes}
\usepackage{stmaryrd}

\usepackage{tikzit}

\tikzstyle{node}=[fill=black, draw=black, shape=circle, scale=0.5]
\tikzstyle{medium_box}=[fill=white, draw=black, shape=rectangle, minimum height=0.8cm, minimum width=0.5cm]

\tikzstyle{pointy}=[->]
\tikzstyle{bluearrow}=[->, fill=none, draw={rgb,255: red,29; green,206; blue,255}, thick]
\tikzstyle{lightnone}=[-, draw={rgb,255: red,191; green,191; blue,191}]

\tikzset{
	Rightarrow/.style={double equal sign distance,>={Implies},->},
	triple/.style={-,preaction={draw,Rightarrow}},
	quadruple/.style={preaction={draw,Rightarrow,shorten >=0pt},shorten >=1pt,-,double,double distance=0.2pt}
}

\newcommand{\Lan}{\operatorname{Lan}}

\newcommand{\profto}{\mathrel{\ooalign{\hfil$\mapstochar$\hfil\cr$\to$\cr}}}

\newcommand{\cat}[1]{\mathcal{#1}}
\newcommand{\ncat}[1]{\mathbf{#1}}
\newcommand{\twocat}[1]{\mathbb{#1}}

\newcommand{\Cat}{\twocat{C}\ncat{at}}

\newcommand{\Set}{\ncat{Set}}

\newcommand{\op}{\mathsf{op}}

\newcommand{\iso}[1][]{\overset{#1}{\cong}}

\font\maljapanese=dmjhira at 2ex
\newcommand{\yo}{\operatorname{\text{\maljapanese\char"48}\!}}

\newcommand{\adj}{\dashv}

\newcommand{\Para}{\ncat{Para}}

\newcommand{\Optic}{\ncat{Optic}}

\newcommand{\action}{\bullet}
\newcommand{\actionn}{\circ}

\newcommand{\comp}{\fatsemi}

\newcommand{\MonCat}{\twocat{M}\ncat{onCat}}

\newcommand{\deloop}{\twocat{B}}

\newcommand{\Actt}{\twocat{A}\ncat{ct}}

\newcommand{\st}{\mathrm{st}}

\newcommand{\equalto}{=\mathrel{\mkern-3mu}=}

\newcommand{\Tamb}{\ncat{Tamb}}
\newcommand{\Prof}{\ncat{Prof}}
\newcommand{\DLens}{\ncat{DLens}}

\newcommand{\Mat}{\twocat{M}\ncat{at}}

\allowdisplaybreaks
\begin{document}
	\begin{abstract}
		Tambara modules are strong profunctors between monoidal categories.
		They've been defined by Tambara in the context of representation theory \cite{tambara2006}, but quickly found their way in applications when it was understood Tambara modules provide a useful encoding of modular data accessors known as mixed optics.
		To suit the needs of these applications, Tambara theory has been extended to profunctors between categories receiving an action of a monoidal category.

		Motivated by the generalization of optics to dependently-typed contexts, we sketch a further extension of the theory of Tambara modules in the setting of actions of double categories (thus doubly indexed categories), by defining them as horizontal natural transformations.
		The theorems and constructions in \cite{pastrostreet2008} relevant to profunctor representation theorem for mixed optics \cite{bryce2020categoricalupdate} are reobtained in this context.
		This reproduces the definition of dependent optics recently put forward by Vertechi \cite{vertechi2022} and Milewski \cite{milewski2022compound}, and hinted at in \cite{braithwaite2021fibre}.
	\end{abstract}

	\maketitle

	\section{Introduction}
In the last two years, and more intensely in the last six months or so, there's been interest in extending the theory of optics to support more strongly typed data structures, giving rise to a many attempts to define so-called `dependent optics'.

The purpose of this search is, on a mathematical level, to unify Myers and Spivak's `$F$-lenses' \cite{spivak2019generalized} with mixed optics \cite{bryce2020categoricalupdate}. The first provide a very flexible and well-motivated framework for dealing with bidirectional systems, further corroborated by their applications in \cite{myers2020double,myers2021book}. The second generalize lenses in many directions, providing a way to speak of bidirectional data accessing in the presence of side-effects and non-cartesian structure, as well as providing a convincing operational semantics for feedback systems \cite{towards_foundations_categorical_cybernetics}.

Some plausible candidates of a common generalization have been defined. The first work the author knows about this has been done by Hedges and Rischel, on what have been called `indexed optics' in \cite{braithwaite2021fibre}. The focus of that definition is on equipping mixed optics with coproducts, a very useful feature that sets apart dependent and simple lenses (and, in that case, everything there is to the generalization: dependent lenses are exactly the coproduct completion of lenses).
In the meantime, Milewski has worked out a very similar definition he called \emph{ommatidia} and then \emph{polynomial optics}---this work can be found in his blog \cite{milewski_polylens}. This essentially coincides with indexed optics, with some minor differences making the latter more general.

In \cite{braithwaite2021fibre} the author and his collaborators expanded on Hedges and Rischel's work on indexed optics in order to find a principled framework for dependent optics.
Two ideas were brought about. First, to define optics `globally' as pullbacks of 2-fibrations of coparameterised morphisms along 2-opfibrations of parameterised morphisms, as opposed to the usual `local' definition which defines objects and hom-sets.
Second, to focus on actions of bicategories. In fact there is a correspondence between actions of bicategories (on categories), (locally discrete) 2-fibrations and 2-indexed categories, which we describe below. Unfortunately, at the time \cite{braithwaite2021fibre} was put out we didn't exactly understood how and if these two ideas could do what we wanted them to do.

In the following months, it became clear that both these intuitions were indeed right.

The `pullback definition' of optics suggests some higher-dimensional structure with a clear operational interpretation is forgot in the usual definition as coends, something described in a particular case by Gavranovi\'c in \cite{gavranovic2022local}.

At the same time, `action of bicategories' have been confirmed as the most promising tool to define dependent optics, since at least two works hingeing on this idea have been published in the meantime.
The first is by Milewski, who defined `compound optics' in \cite{milewski2022compound} to be a vast and principled generalization of his ommatidia and copresheaf optics.
The second is the very recent work by Vertechi \cite{vertechi2022}, which takes seriously the idea of actions of bicategories to define `dependent optics'. He proposes various new definitions of dependently-typed modular data accessors generalizing dependent lenses and works out a tentative Tambara representation theorem.

In this note, which we developed independently from the works just mentioned, we extend the theory of Tambara modules in the setting of actions of double categories (also known as doubly indexed categories), by defining them as horizontal natural transformations between doubly indexed categories, a remarkably simple definition.
The theorems and constructions in \cite{pastrostreet2008} relevant to profunctor representation theorem for mixed optics \cite{bryce2020categoricalupdate} are reobtained in this context.
If the abstract definitions were independently discovered, we lacked confidence in their practical usefulness until Milewski's published explicit calculations for indexed optics (aka polylenses, aka  ommatidia).
We have been convinced that the definition we had was indeed promising only when Vertechi shown that it does indeed capture dependent lenses in the expected way (as we lacked the idea used in~\cref{ex:dlens} and first appeared in \cite{vertechi2022}).

	\section{From actions of monoids to actions of bicategories}
The action of a monoid $M$ on a set $X$ is a well-known structure: it's a map that associates an endomorphism $m \action - : X \to X$ to each $m \in M$.
The action of a monoid might as well be presented as a \emph{representation of $M$ in $\Set$}, i.e.~a functor
\begin{equation}
	\deloop M \longto \Set
\end{equation}
where $M$ has been correctly identified as a one-object category, whose only object hits the set acted upon by $M$ and whose endomorphisms hit endomorphisms of $X$, as desired.
Therefore it is easy to generalize this to categories (i.e.~`monoids with many objects'): one simply replaces $\deloop M$ with any other category $\cat C$.
The resulting functor $\cat C \to \Set$ is more commonly known as a \emph{copresheaf over $\cat C$}, or \emph{$\cat C$-indexed set}, but can be rigthfully considered to be an \emph{action} or \emph{representation of $\cat C$}.

This point of view might feel a bit artificial.
Actions of a monoid are quite humble structures, definable in any monoidal category, whereas here we invoke the entire category of sets to define one.
The uneasiness we feel is that of having moved from an \emph{internal} definition (`an action is a map of sets such that...') to an \emph{external} definition (`an action is a functor to sets').
To ease our discomfort, we seek an `internal' presentation of actions of categories.

In what follows, let $\cat E$ be a sufficiently structured category. We are mostly going to be interested in $\cat E = \Set$, but at least for the following definition any finitely complete category (i.e.~admitting finite limits) suffices.

\begin{definition}
\label{def:action-of-category}
	Let $\cat M : M_1 \rightrightarrows M_0$ be an internal category\footnotemark\ in $\cat E$. Let $X$ be an object over $M_0$, as witnessed by a map $\pi : X \to M_0$.
	An \emph{action of $\cat M$ on $X$} is a map:
	\begin{equation}
		{\action} : (o, o' : M_0) \to (m : M_1\,o\,o') \to X\, o \to X\, o'
	\end{equation}
	satisfying a form of `functoriality' (here $o \nto{m} o' \nto{n} o''$ are consecutive morphisms in $\cat M$):
	\begin{equation}
		1_{o} \action x = x, \qquad m \action (n \action x) = (m \comp n) \action x.
	\end{equation}
\end{definition}
\footnotetext{
	For reference purposes, we mention the definition of internal category:
	\begin{definition}
	\label{def:internal-category}
		An \emph{internal category} $\cat M$ is given by (1) an object $M_0$ of `objects',
		(2) an object $M_1$ of `arrows',
		(3) two maps $s,t:M_1 \to M_0$ assigning a source and target to each arrow,
		(4) a map $1: M_0 \to M_1$ picking the identity arrows and
		(5) a map $\comp : M_1 \times_{M_0} M_1 \to M_1$ composing consecutive pairs of arrows (the pullback is taken over the cospan $M_1 \nto{s} M_0 \nfrom{t} M_1$);
		satisfying the usual associativity and unitality of composition in a category.
	\end{definition}
}

Equivalently, $\action$ is an arrow like the blue one in this diagram:
\begin{equation}
	\label{diag:action-of-category}
	\begin{tikzcd}
		\textcolor{rgb,255:red,36;green,56;blue,204}{X} & \textcolor{rgb,255:red,36;green,56;blue,204}{M \times_{M_0} X} & X \\
		{M_0} & M & {M_0}
		\arrow["s"', from=2-2, to=2-3]
		\arrow["\pi", from=1-3, to=2-3]
		\arrow["{p_1}", from=1-2, to=1-3]
		\arrow["{p_2}", from=1-2, to=2-2]
		\arrow["\lrcorner"{anchor=center, pos=0.125}, draw=none, from=1-2, to=2-3]
		\arrow["t", from=2-2, to=2-1]
		\arrow["\action"', color={rgb,255:red,36;green,56;blue,204}, from=1-2, to=1-1]
		\arrow["\pi"', from=1-1, to=2-1]
	\end{tikzcd}
\end{equation}
Notice the right square is a pullback but the left square isn't.

If the map $\pi$ is `nice enough' (or the category $\cat E$ is nice enough), $X$ can be considered a dependent type
\begin{equation*}
	o : M_0 \vdash X\, o\ \text{type}.
\end{equation*}
In $\Set$, this means $X$ is actually an ${M_0}$-indexed family of sets, woven together by the action of $\cat M$:
\begin{figure}[H]
	\includegraphics[width=.4\textwidth]{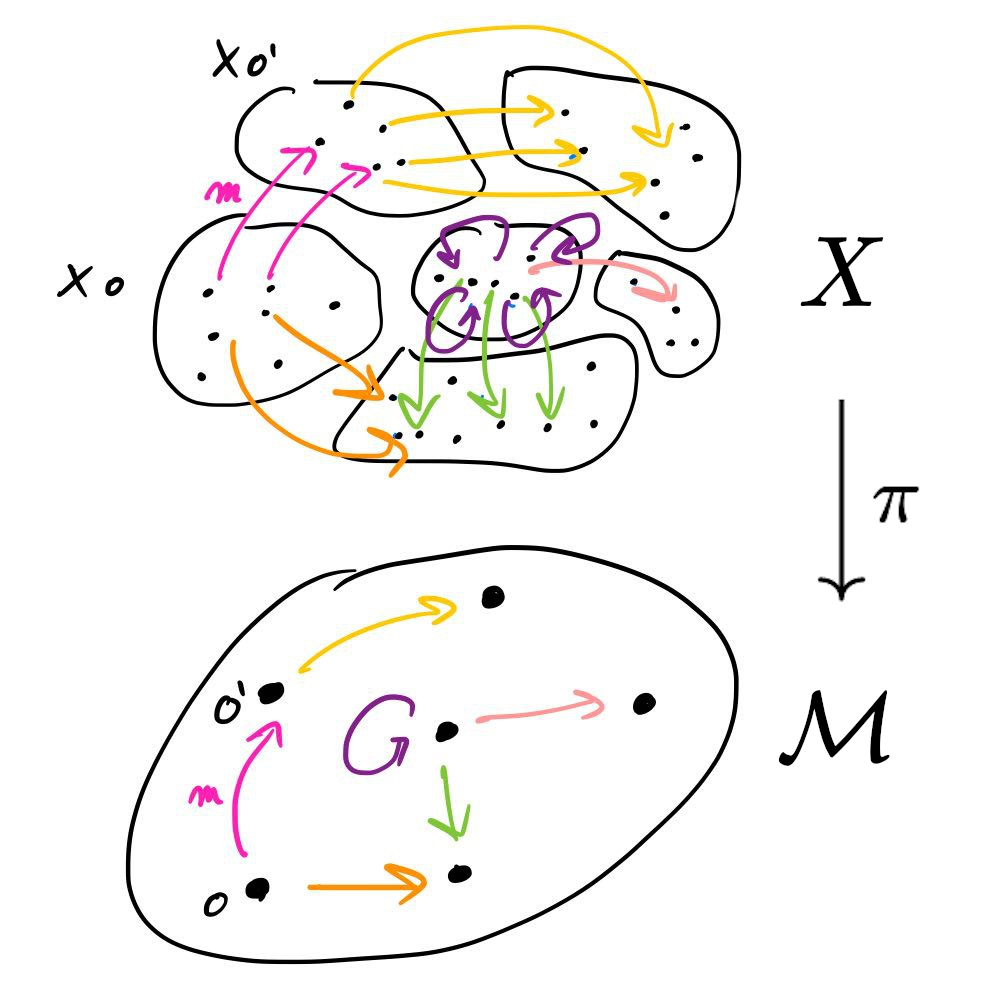}
\end{figure}
Hence in the case $\cat E = \Set$, where an internal category is the same thing as a small category, one readily pass from this definition to that of a copresheaf over $\cat M$:
\begin{equation}\begin{tikzcd}
	{\cat M} & \Set \\[-6ex]
	o & {X\, o} \\
	{o'} & {X\, o'}
	\arrow[from=1-1, to=1-2]
	\arrow[""{name=0, anchor=center, inner sep=0}, "m"', from=2-1, to=3-1]
	\arrow[""{name=1, anchor=center, inner sep=0}, "{m \action -}", from=2-2, to=3-2]
	\arrow[shorten <=10pt, shorten >=10pt, maps to, from=0, to=1]
\end{tikzcd}\end{equation}

So the definition of action of a category \emph{can} be given internally. There is another intermediate way one can express an action, which is to define it as an \emph{opfibration} over $\cat M$.
In fact, by Grothendieck construction, the data of $\action : \cat M \to \Set$ can be bundled into a single functor with the structure of a (discrete) opfibration:
\begin{equation}
\begin{tikzcd}
	{\Para_\action(X)} \\[-3ex]
	{\cat M}
	\arrow["p", from=1-1, to=2-1]
\end{tikzcd}
\end{equation}
The total category of this fibration is given by pairs $(o : \cat M, x : X\,o)$, hence by elements of $X$, with morphisms $(m, =) : (o,x) \to (o',x')$ given by maps $m: o \to o'$ such that $m \action x = x'$.
Sometimes $\Para_\action(X)$ is called the `action category' associated to $\action$.
The lifting structure on the projection functor $p$, making it into an opfibration, is given by
\begin{equation}
\begin{tikzcd}
	x & {m \action x} \\[-3ex]
	o & {o'}
	\arrow["m", from=2-1, to=2-2]
	\arrow["{(m, 1_{m \action x})}", from=1-1, to=1-2]
	\arrow[dotted, maps to, from=1-1, to=2-1]
	\arrow[dotted, maps to, from=1-2, to=2-2]
\end{tikzcd}
\end{equation}
Hence `actions' of $\cat M$ seem to admit three equivalent presentations:
\begin{enumerate}
	\item as copresheaves $\cat M \to \Set$,
	\item as (discrete) opfibrations $\cat P \to \cat M$,
	\item as pairs of maps $X \to M_0$, $M_1 \times_{M_0} X \to X$ satisfying some laws.
\end{enumerate}

	Notice the adjective \emph{discrete} for the opfibration description of an action. This discreteness comes form the discreteness of $X$: it's just a set (or an object of $\cat E$). One can go beyond this and ask for $X$ to be itself a category (or an internal category in $\cat E$).
	Then, externally, the action of $\cat M$ would correspond to a functor $\cat M \to \Cat$ or an unqualified opfibration $\cat P \to \cat M$.

% ...

By further increasing the dimensionality of the objects involved, we move from actions of monoidal categories on categories (also known as \emph{actegories}), to actions of \emph{bicategories} on categories:

\begin{definition}
	The \textbf{action of a bicategory} $\twocat M$ on a category $\cat X$ is, equivalently:
	\begin{enumerate}
		\item an indexed category $\twocat M \to \Cat$,
		\item a (locally discrete) 2-opfibration $\twocat E \to \twocat M$,
		\item a pair of functors $\cat X \to \cat M_0$ and $\cat M_1 \times_{\cat M_0} \cat X \to \cat X$ equipped with further strucure morphisms categorifying the laws in~\cref{def:action-of-category} (see \cite{bakovic2009simplicial}).
	\end{enumerate}
\end{definition}

\subsection{The canonical action of spans}
\label{subsec:spans}
The `ur-example' of action of a bicategory is the action of $\Span(\cat I)$ `on $\cat I$' for a category with finite limits $\cat I$.
Remember $\Span(\cat I)$ is a bicategory whose objects are the same as $\cat I$, while morphisms $i,j:\cat I$ are given by spans
\begin{equation}
\begin{tikzcd}[sep=4ex]
	& z \\[-3ex]
	i && j
	\arrow[from=1-2, to=2-1]
	\arrow[from=1-2, to=2-3]
\end{tikzcd}
\end{equation}
composing by pullback.
The 2-morphisms are given by apex-wise morphisms commuting with the legs:
\begin{equation}
\begin{tikzcd}[sep=5ex]
	& z \\[-3ex]
	i && j \\[-3ex]
	& w
	\arrow[from=1-2, to=2-1]
	\arrow[from=1-2, to=2-3]
	\arrow[from=3-2, to=2-1]
	\arrow[from=3-2, to=2-3]
	\arrow[Rightarrow, from=1-2, to=3-2, "\alpha"]
\end{tikzcd}
\end{equation}
The `canonical' action of $\Span(\cat I)$ is given by the following functor:
\begin{equation}
\label{eq:canonical-span-action}
\begin{tikzcd}
	&[-10ex] {\Span(\cat I)} &[-10ex]&& \Cat \\[-6ex]
	& i &&& {\cat I/i} \\
	z && w \\
	& j &&& {\cat I/j}
	\arrow["q"', from=3-1, to=4-2]
	\arrow["p", from=3-1, to=2-2]
	\arrow["s", from=3-3, to=4-2]
	\arrow["r"', from=3-3, to=2-2]
	\arrow["\alpha"{description}, Rightarrow, from=3-1, to=3-3]
	\arrow[""{name=0, anchor=center, inner sep=0}, "{q_* p^*}"{description}, curve={height=24pt}, from=2-5, to=4-5]
	\arrow[""{name=1, anchor=center, inner sep=0}, "{s_*r^*}"{description}, curve={height=-24pt}, from=2-5, to=4-5]
	\arrow["L", from=1-2, to=1-5]
	\arrow["{\alpha^*}", shorten <=10pt, shorten >=10pt, Rightarrow, from=0, to=1]
	\arrow[shorten <=9pt, shorten >=13pt, maps to, from=3-3, to=0]
\end{tikzcd}
\end{equation}
This is also known as the action of \emph{linear polynomial functors in many indeterminates} \cite{gambino2013polynomial}, since it's given by associating each span with the functor that pulls back and then push forward along the two legs of the span.
\begin{equation}\begin{tikzcd}
	& {\cat I/z} \\[-3.5ex]
	{\cat I/i} & z & {\cat I/j} \\[-3.5ex]
	i && j
	\arrow["{p^*}", from=2-1, to=1-2]
	\arrow["{q_*}", from=1-2, to=2-3]
	\arrow["p"', from=2-2, to=3-1]
	\arrow["q", from=2-2, to=3-3]
	\arrow[dotted, maps to, from=2-1, to=3-1]
	\arrow[dotted, maps to, from=1-2, to=2-2]
	\arrow[dotted, maps to, from=2-3, to=3-3]
\end{tikzcd}\end{equation}
The reason these are called linear is that pushforward is a kind of `fibrewise sum', hence the notation $\sum_q$ sometimes used for $q_*$ (for instance in \cite{schreiber2014quantization}).

\begin{remark}
	The action of $\Span(\cat I)$ is secretely given by the fibrewise self-action of the codomain fibration of $\cat I$, which is a fancy way to say $L$ is packaging all the actions
	\begin{equation}
		\times_i : \cat I/i \times \cat I/i \longto \cat I/i
	\end{equation}
	given by cartesian product in $\cat I$.
	Since pullbacks paste, reindexing are monoidal between these categories.
	We conclude that the indexed category $\cat I/- : \cat I^\op \to \Cat$ is actually an \emph{indexed monoidal category} \cite{moeller_vasilakopoulou_monoidal_grothendieck} with left adjoints to all reindexing functors.
	Furthermore, one may notice we can weaken even more the structure: instead of a full-fledged monoidal category, we could ask to have an actegory $(\cat M_i, \cat C_i, \action_i)$ for each $i: \cat I$, with linear reindexing with left adjoints. This will be the subject of~\cref{ex:fibre-optics}
\end{remark}

\begin{remark}
	We emphasize that the natural object of the action of a bicategory is not simply a category, but a category \emph{over the set of objects of the acting bicategory}.
	Such a dependency is very often decorated by more structure. For example, $L$ above can be considered to be $\Span(\cat I)$ \emph{acting on the codomain fibration of $\cat I$}, hence the dependency on the objects of $\cat I$ is actually a whole fibration on the entire $\cat I$.
\end{remark}

\subsection{The action of profunctors}
\label{subsec:profs}
In \cite{milewski2022compound}, Milewski defines an action analogous to the one of $\Span(\cat I)$ to define `compound optics'. The action is that of $\Prof$ on $\Cat$, or better, on the 2-fibration of discrete fibrations over small categories (compare this with: the fibration of morphisms over objects of $\cat I$).
Explicitly, $\Prof$ is taken to act on the 2-indexed category $[-,\Set] : \Cat^\op \to \Cat$ that sends a category $\cat C$ to its category of copresheaves.
\begin{diagram}
	\Prof && \Cat \\[-5ex]
	{\cat C} && {[\cat C, \Set]} \\
	{\cat D} && {[\cat D, \Set]}
	\arrow[""{name=0, anchor=center, inner sep=0}, "Q","\shortmid"{marking}, curve={height=-18pt}, from=2-1, to=3-1]
	\arrow[""{name=1, anchor=center, inner sep=0}, "P"', "\shortmid"{marking}, curve={height=18pt}, from=2-1, to=3-1]
	\arrow[""{name=2, anchor=center, inner sep=0}, "{P\bullet }"', curve={height=18pt}, from=2-3, to=3-3]
	\arrow[""{name=3, anchor=center, inner sep=0}, "Q\bullet", curve={height=-18pt}, from=2-3, to=3-3]
	\arrow["\bullet", from=1-1, to=1-3]
	\arrow["\alpha\bullet", shift right=1, shorten <=5pt, shorten >=5pt, Rightarrow, from=2, to=3]
	\arrow["\alpha", shift right=1, shorten <=5pt, shorten >=5pt, Rightarrow, from=1, to=0]
	\arrow[shorten <=28pt, shorten >=28pt, maps to, from=0, to=2]
\end{diagram}
The functor $P\action$ is defined by what Milewski suggests is a `categorified vector-matrix product':
\begin{equation}
\label{eq:compound-action}
	P \action a = \int^{c: \cat C} a(c) \times P(c, -)
\end{equation}
The action on a natural transformation $\alpha$ is defined by functoriality of the coend.

	\section{Tambara modules for actions of double categories}
\label{sec:tambara}
For actions of monoids, there is a notion of `lax equivariant relation', that is, if $X$ and $Y$ are $M$-sets, a relation $P : X \profto Y$ such that
\begin{equation}
	\forall m \in M, \qquad P(x, y) \word{implies} P(m \action x, m \action y).
\end{equation}
This notion famously categorifies to Tambara modules \cite{tambara2006}, which are now \emph{profunctors} $P : \cat X \profto \cat Y$, for $\cat X$ and $\cat Y$ left $\cat M$-actegories, equipped with a family of morphisms:
\begin{equation}
	\st_{m,x,y} : P(x,y) \longto P(m \action x, m \action y)
\end{equation}
which are dinatural in $m$ and natural in $x$ and $y$.

How do we generalize this to actions of bicategories?
Let us think first of the 0-dimensional case. Suppose now $X$ and $Y$ are sets over $M_0$, the set of objects of a category $\cat M : M_1 \rightrightarrows M_0$. Remember this is equivalent to have $M_0$-indexed families of sets `woven together' by morphisms of $\cat M$. Being $X$ and $Y$ naturally defined over $M_0$, it makes sense to also define $P$ as a `relation over $M_0$', that is, a family of relations:
\begin{equation}
	P: (o : M_0) \to (X\, o \profto Y\, o)
\end{equation}
Of course the delicate arrangement of $P$ (elements over $o$ are related only to elements over $o$) can be forgotten, and one can consider $P$ as a relation $X \profto Y$. This, however, would make impossible to define the compatibility of $P$ with the $\cat M$-actions on $X$ and $Y$:
\begin{equation}
	\forall o \nto{f} o' : \cat M, \qquad (P\, o)(x, y) \word{implies} (P\, o')(m \action x, m \action y).
\end{equation}
Hence $P$ is crucially a relation \emph{over $M_0$} `woven together' by morphisms of $\cat M$, in perfect analogy with the notion of action of $\cat M$ itself.

Categorifying this makes us summon a $\cat M$-indexed family of profunctors, which amounts to a profunctor
\begin{equation}
	P_o : \cat X(o) \profto \cat Y(o)
\end{equation}
for each $o : \cat M$ and a natural transformation
\begin{equation}
	\st_f : P_o \longtwoto P_{o'}(m \action -, m \action -)
\end{equation}
for each $f: o \to o'$.

This could be made directly into a definition (as done in \cite{vertechi2022}), but we'd be missing out in elegance.
Indeed, the peculiar form of $\st_f$ should make anyone familiar with the proarrow equipment of categories and profunctors suspicious.
Such suspicion is aggravated by the fact that already the `nuts and bolts' definition of action of a bicategory seems to beg for the bicategory $\cat M$ to be promoted to a double category.
In fact the analogoue definition one dimension lower presents the acting category as an internal category in $\Set$, so it feels much more natural to promote `internal category in $\Set$' (i.e.~`small category') to `internal category in $\Cat$' (i.e.~double category).
Everything conjures, therefore, to the following slightly more general definitions:

\begin{definition}
	The \textbf{action of a double category} $\twocat M$ on a category $\cat X$ is, equivalently:
	\begin{enumerate}
		\item a doubly indexed category $\twocat M \to \Cat$,
		\item a (locally discrete) double 2-opfibration $\twocat E \to \twocat M$,
		\item a pair of functors $\cat X \to \cat M_0$, $\cat M_1 \times_{\cat M_0} \cat X \to \cat X$ equipped with further structure morphisms categorifying~\cref{def:action-of-category} (see \cite{moskaliuk1998double}).
	\end{enumerate}
\end{definition}

\begin{definition}
\label{def:tambara-module}
	Let $\twocat M$ be a double category.
	Suppose $\cat X$ and $\cat Y$ are doubly indexed categories $\twocat M \to \Cat$.
	Then a \textbf{Tambara module} $(P, \st)$ is a horizontal natural transformation $\cat X \twoto \cat Y$:
	\begin{diagram}[sep=4ex]
		{\cat X(o)} & {\cat Y(o)} \\
		{\cat X(o')} & {\cat Y(o')}
		\arrow["{P_o}", "\shortmid"{marking}, from=1-1, to=1-2]
		\arrow["{P_{o'}}"', "\shortmid"{marking}, from=2-1, to=2-2]
		\arrow["{m \action -}"', from=1-1, to=2-1]
		\arrow["{m \action -}", from=1-2, to=2-2]
		\arrow["\st"', shorten <=7pt, shorten >=7pt, Rightarrow, from=1-2, to=2-1]
	\end{diagram}
\end{definition}

The naturality of $(P, \st)$ amounts to ask, among the other things (see, \emph{mutatis mutandis},~\cite[Definition 1.6]{pare2011yoneda}, which defines \emph{vertical} natural transformations):
\begin{diagram}[sep=4ex]
	\label{diag:tambara-law1}
	{\cat X(o)} & {\cat Y(o)} && {\cat X(o)} & {\cat Y(o)} \\
	{\cat X(o)} & {\cat Y(o)} && {\cat X(o)} & {\cat Y(o)}
	\arrow["{P_o}", "\shortmid"{marking}, from=1-1, to=1-2]
	\arrow["{P_o}"', "\shortmid"{marking}, from=2-1, to=2-2]
	\arrow["{1_o \action -}"', from=1-1, to=2-1]
	\arrow[""{name=0, anchor=center, inner sep=0}, "{1_o \action -}", from=1-2, to=2-2]
	\arrow["\st"', shorten <=7pt, shorten >=7pt, Rightarrow, from=1-2, to=2-1]
	\arrow[""{name=1, anchor=center, inner sep=0}, Rightarrow, no head, from=1-4, to=2-4]
	\arrow[Rightarrow, no head, from=1-5, to=2-5]
	\arrow["{P_o}", "\shortmid"{marking}, from=1-4, to=1-5]
	\arrow["{P_o}"', "\shortmid"{marking}, from=2-4, to=2-5]
	\arrow["{1_{P_o}}"', shorten <=7pt, shorten >=7pt, Rightarrow, no head, from=1-5, to=2-4]
	\arrow[shorten <=30pt, shorten >=24pt, Rightarrow, no head, from=0, to=1]
\end{diagram}
\begin{diagram}[sep=4ex]
\label{diag:tambara-law2}
	{\cat X(o)} & {\cat Y(o)} &&& {\cat X(o)} & {\cat Y(o)} \\
	{\cat X(o')} & {\cat Y(o')} \\
	{\cat X(o'')} & {\cat Y(o'')} &&& {\cat X(o'')} & {\cat Y(o'')}
	\arrow["{P_o}", "\shortmid"{marking}, from=1-1, to=1-2]
	\arrow["{P_{o'}}"', "\shortmid"{marking}, from=2-1, to=2-2]
	\arrow["{m \action -}"', from=1-1, to=2-1]
	\arrow["{m \action -}", from=1-2, to=2-2]
	\arrow["\st"', shorten <=8pt, shorten >=8pt, Rightarrow, from=1-2, to=2-1]
	\arrow["{n \action -}"', from=2-1, to=3-1]
	\arrow["{n \action -}", from=2-2, to=3-2]
	\arrow["{P_{o''}}"', "\shortmid"{marking}, from=3-1, to=3-2]
	\arrow["{P_o}", "\shortmid"{marking}, from=1-5, to=1-6]
	\arrow["{P_{o''}}"', "\shortmid"{marking}, from=3-5, to=3-6]
	\arrow[""{name=0, anchor=center, inner sep=0}, "{(m \comp n) \action -}"', from=1-5, to=3-5]
	\arrow["{(m \comp n) \action -}", from=1-6, to=3-6]
	\arrow["\st"', shorten <=8pt, shorten >=8pt, Rightarrow, from=2-2, to=3-1]
	\arrow["\st"', shorten <=14pt, shorten >=14pt, Rightarrow, from=1-6, to=3-5]
	\arrow[shorten <=12pt, shorten >=47pt, Rightarrow, no head, from=2-2, to=0]
\end{diagram}

\begin{remark}
	Just for fun, the transpose definition would be that of a family of functors $F_o: \cat X(o) \to \cat Y(o)$ that commute with the horizontal part of the action of $\twocat M$, yielding squares
	\begin{diagram}
		{\cat X(o)} & {\cat X(o')} \\
		{\cat Y(o)} & {\cat Y(o')}
		\arrow["{- \actionn m \action =}", "\shortmid"{marking}, from=1-1, to=1-2]
		\arrow["{- \actionn m \action =}"', "\shortmid"{marking}, from=2-1, to=2-2]
		\arrow["{F_o}"', from=1-1, to=2-1]
		\arrow["{F_{o'}}", from=1-2, to=2-2]
		\arrow["\alpha"', shorten <=7pt, shorten >=7pt, Rightarrow, from=1-2, to=2-1]
	\end{diagram}
	given again by natural transformations
	\begin{equation}
		\alpha_m : x \actionn m \action y \longto F_o(x) \actionn m \action F_{o'}(y)
	\end{equation}
	whose interpretation is obscure to the author.
\end{remark}

In the case $\twocat M$ is not only a `discrete' double category, but a \emph{one-object} discrete double category, then this recovers exactly `classical' Tambara modules between actions of monoidal categories \cite{bryce2020categoricalupdate}.
All squares look like this (notice the lack of indexing over objects):
\begin{diagram}
	{\cat X} & {\cat Y} \\
	{\cat X} & {\cat Y}
	\arrow["P", "\shortmid"{marking}, from=1-1, to=1-2]
	\arrow["P"', "\shortmid"{marking}, from=2-1, to=2-2]
	\arrow["{m \action -}"', from=1-1, to=2-1]
	\arrow["{m \action -}", from=1-2, to=2-2]
	\arrow["st"', shorten <=6pt, shorten >=6pt, Rightarrow, from=1-2, to=2-1]
\end{diagram}
and the laws~\eqref{diag:tambara-law1}--\eqref{diag:tambara-law2} encode the usual laws for a Tambara structure.
If $\twocat M$ is a bicategory encoded as an horizontally discrete double category, then nothing much changes, and laws~\eqref{diag:tambara-law1}--\eqref{diag:tambara-law2} are the only ones which remain non-trivial.

\subsection{Pastro--Street theory.}
In \cite{pastrostreet2008}, a free-forgetful-cofree adjunction is constructed between profunctors $P:\cat X \profto \cat Y$ and Tambara modules $(P, \st): (\cat X, \action) \profto (\cat Y, \action)$.
To reconstruct this adjunction, we need first of all to define what's a `profunctor' without the Tambara structure.

We expect it to be the same as our initial ansatz, i.e.~simply a profunctor over $\cat M_0$.
To check this, we ask given two doubly indexed categories $\cat X, \cat Y: \twocat M \to \Cat$, what are the corresponding `bare' categories to be used as domain and codomain for such a profunctor?
It's all that left if we `forget' all the data in $\cat X$ associated with the action of non-trivial morphisms in $\twocat M$, i.e.~the (`singly') indexed category
\begin{eqalign}
	&\cat X_0 := \cat M_0 \into \twocat M \nlongto{\cat X} \Cat,\\
	&\cat Y_0 := \cat M_0 \into \twocat M \nlongto{\cat Y} \Cat
\end{eqalign}
One quickly realizes a `Tambara module' for this trivialization is indeed a profunctor over $\cat M_0$ (the strength is forced to be the identity by~\eqref{diag:tambara-law1}):
\begin{diagram}[sep=4ex]
	{\cat X(o)} & {\cat Y(o)} \\
	{\cat X(o)} & {\cat Y(o)}
	\arrow["{P_o}", "\shortmid"{marking}, from=1-1, to=1-2]
	\arrow["{P_o}"', "\shortmid"{marking}, from=2-1, to=2-2]
	\arrow[Rightarrow, no head, from=1-1, to=2-1]
	\arrow[Rightarrow, no head, from=1-2, to=2-2]
\end{diagram}

Given two $\twocat M$-actions $\cat X$, $\cat Y$, we have thus two categories:
\begin{enumerate}
	\item The category $\Tamb(\cat X, \cat Y)$ of Tambara modules in the sense of definition~\cref{def:tambara-module}, having as morphisms (natural) families of (natural) transformations $\beta_o : P_o \longtwoto Q_o$ such that all equations
	\begin{diagram}[sep=4ex]
		{\cat X(o)} & {\cat Y(o)} && {\cat X(o)} & {\cat Y(o)} \\
		\\
		{\cat X(o')} & {\cat Y(o')} && {\cat X(o')} & {\cat Y(o')}
		\arrow[""{name=0, anchor=center, inner sep=0}, "{Q_o}"', "\shortmid"{marking}, curve={height=12pt}, from=1-1, to=1-2]
		\arrow["{Q_{o'}}"', "\shortmid"{marking}, from=3-1, to=3-2]
		\arrow["{m \action -}"', from=1-1, to=3-1]
		\arrow[""{name=1, anchor=center, inner sep=0}, "{m \action -}", from=1-2, to=3-2]
		\arrow["{\st_Q}", shift left=2, shorten <=14pt, shorten >=14pt, Rightarrow, from=1-2, to=3-1]
		\arrow[""{name=2, anchor=center, inner sep=0}, "{P_o}", "\shortmid"{marking}, curve={height=-12pt}, from=1-1, to=1-2]
		\arrow[""{name=3, anchor=center, inner sep=0}, "{m \action -}"', from=1-4, to=3-4]
		\arrow["{m \action-}", from=1-5, to=3-5]
		\arrow["{P_o}", "\shortmid"{marking}, from=1-4, to=1-5]
		\arrow[""{name=4, anchor=center, inner sep=0}, "{P_{o'}}", "\shortmid"{marking}, curve={height=-12pt}, from=3-4, to=3-5]
		\arrow[""{name=5, anchor=center, inner sep=0}, "{Q_{o'}}"', "\shortmid"{marking}, curve={height=12pt}, from=3-4, to=3-5]
		\arrow["{\st_P}"', shift right=2, shorten <=14pt, shorten >=14pt, Rightarrow, from=1-5, to=3-4]
		\arrow["{\beta_o}", shorten <=3pt, shorten >=3pt, Rightarrow, from=2, to=0]
		\arrow["{\beta_{o'}}", shorten <=3pt, shorten >=3pt, Rightarrow, from=4, to=5]
		\arrow[shorten <=32pt, shorten >=32pt, Rightarrow, no head, from=1, to=3]
	\end{diagram}
	hold.
	\item The category $\Prof(\cat X_0, \cat Y_0)$ of profunctors over $\cat M_0$, having as morphisms (natural) families of (natural) transformations $\beta_o : P_o \longtwoto Q_o$.
\end{enumerate}

We are looking for a triple of adjoint functors
\begin{diagram}
	{\Tamb(\cat X, \cat Y)} && {\Prof(\cat X_0, \cat Y_0)}
	\arrow[""{name=0, anchor=center, inner sep=0}, "U"{description}, from=1-1, to=1-3]
	\arrow[""{name=1, anchor=center, inner sep=0}, "\Theta"{description}, shift right=5, curve={height=6pt}, from=1-3, to=1-1]
	\arrow[""{name=2, anchor=center, inner sep=0}, "\Psi"{description}, shift left=5, curve={height=-6pt}, from=1-3, to=1-1]
	\arrow["\dashv"{anchor=center, rotate=-90}, draw=none, from=1, to=0]
	\arrow["\dashv"{anchor=center, rotate=-90}, draw=none, from=0, to=2]
\end{diagram}
where $U$ forgets the structure ($\st$) as well as the properties of a Tambara module.

One can prove this directly, by using the same constructions expounded in \cite[Chapter 5]{roman2020profunctor} and derived from Pastro and Street's original proof.
That path involves quite a lot of coend calculus and the strategy is to prove that $U\Theta$ is a comonad on $\Prof(\cat X_0, \cat Y_0)$ and $U\Psi$ is a left adjoint monad.
Alternatively, notice $U$ is really pullback along $\cat M_0 \into \twocat M$ and thus $\Theta$ and $\Psi$ are given by left and right Kan extension.


\begin{proposition}
	The comonad $U\Theta$ is defined on a given $P : \Prof(\cat X_0, \cat Y_0)$ as follows:
	\begin{equation}
		(U\Theta P)_o(x, y) = \int_{o' : \cat M_0} \int_{m : \twocat M(o,o')} P_{o'}(m \action x, m \action y)
	\end{equation}
	for each $o:\cat M_0$.
\end{proposition}

\begin{proposition}
	The right adjoint $U\Psi$ to $U\Theta$ is defined on a given $Q : \Prof(\cat X_0, \cat Y_0)$ as follows:
	\begin{equation}
		(U\Psi Q)_{o'}(x, y) = \int^{o : \cat M_0} \int^{m : \twocat M(o,o')} \int^{a : \cat X(o), b : \cat Y(o)} \cat X(o')(m \action a, x) \times Q_o(a,b) \times \cat Y(o')(y, m \action b)
	\end{equation}
	for each $o' : \cat M_0$.
\end{proposition}

Once $\Psi$ and $\Theta$ are defined, the following is an easy corollary:

\begin{theorem}[Dependent profunctor representation theorem]
\label{th:tambara-rep}
	Let $o,o' : \twocat M$, let $a : \cat X(o)$, $b: \cat Y(o)$, $s:\cat X(o')$, $t:\cat Y(o')$
	\begin{equation}
	\label{eq:profunctor-rep}
		\int_{P : \Tamb(\cat X, \cat Y)} \Set(P_o(a,b), P_{o'}(s,t))
		\iso
		\int^{m : \twocat M(o,o')} \cat X(o')(s, m \action a) \times \cat Y(o')(m \action b, t)
	\end{equation}
\end{theorem}

In the theory of optics, the sets appearing in~\eqref{eq:profunctor-rep} gather as hom-sets of a category, the category of optics.
Thus it feels warranted to do the same here, calling \textbf{dependent optics} the category $\Optic_{\cat X, \cat Y}$ whose objects are given by triples $(o : \cat M, a : \cat X(o), b : \cat Y(o))$ and whose hom-sets are given by the sets appearing in the above isomorphism.

This exact same definition appears in \cite{vertechi2022}, where~\cref{th:tambara-rep} is proved in §4.2, albeit in a different form.

	\section{Examples}
\label{sec:examples}
Arguably the most important example of dependent optics (or better: their most illustrious progenitor) are dependent lenses, which, as we are going to prove, arise from the canonical action of (the bicategory) $\Span(\cat I)$ on the self-indexing of $\cat I$ described in~\cref{subsec:spans}.
This has been worked out by Vertechi in \cite[§3.1]{vertechi2022}. The key observation is the following lemma (unstated in \emph{ibid.}):

\begin{lemma}
\label{lemma:dlens}
	Let $\cat I$ be a finitely complete category.
	Let $v : s \to j$, $u : a \to i$ be maps in $\cat I$, and let $m : i \nfrom{p} z \nto{q} j$ be a span in $\cat I$.
	Then
	\begin{equation}
		\cat I/j(v, q_* p^* u) \iso \sum_{f \in \cat I(s, a)} \cat I/(i \times j)(\langle f \comp u, v \rangle, \langle p, q \rangle)
	\end{equation}
\end{lemma}
\begin{proof}
	Observe an arrow $\phi \in \cat I/j(v, q_* p^* u)$ is given by an arrow $\phi :s \to p^*a$ fitting in the following commutative diagram:
	\begin{diagram}
		a & {p^*a} && s \\
		i & z & j
		\arrow["q"', from=2-2, to=2-3]
		\arrow["p", from=2-2, to=2-1]
		\arrow["u"', from=1-1, to=2-1]
		\arrow[from=1-2, to=1-1]
		\arrow[from=1-2, to=2-2]
		\arrow["\lrcorner"{anchor=center, pos=0.125, rotate=-90}, draw=none, from=1-2, to=2-1]
		\arrow["{q_*p^*u}"{description}, from=1-2, to=2-3]
		\arrow["v"{description}, from=1-4, to=2-3]
		\arrow["\phi"', dotted, from=1-4, to=1-2]
	\end{diagram}
	Since $\phi$ is a map towards the pullback of $u$ and $p$ over $z$, it contains the same data as a map $f:s \to a$ and a map $g:s \to z$ such that $g \comp p = f \comp u$. But we also want the triangle on the right to commute, hence the map $g$ has to satisfy $g \comp q = v$ too.
	We thus draw another diagram, using the observation that the span $m : i \nfrom{p} z \nto{q} j$ is equivalently given as $\langle p, q \rangle : z \to i \times j$:
	\begin{diagram}[sep=5ex]
		&[-2ex]& s \\
		&& z \\
		a && {i \times j} \\
		& i && j
		\arrow["f"', curve={height=12pt}, from=1-3, to=3-1]
		\arrow["u"', from=3-1, to=4-2]
		\arrow["{\pi_i}", from=3-3, to=4-2]
		\arrow["{\pi_j}"', from=3-3, to=4-4]
		\arrow["v", curve={height=-24pt}, from=1-3, to=4-4]
		\arrow["{\langle p, q \rangle}", from=2-3, to=3-3]
		\arrow["{\exists! \langle f \comp u, v \rangle}"'{pos=0.7}, curve={height=25pt}, dashed, from=1-3, to=3-3]
		\arrow["g", dotted, from=1-3, to=2-3]
	\end{diagram}
	The maps $f$ and $v$ can be packaged into a unique map $\langle f \comp u, v \rangle : s \to i \times j$, by universal property of $i \times j$. It's now easy to see a map $\langle f\comp u, v \rangle \to \langle p,q\rangle$ in $\cat I/i \times j$ encodes exactly the data and the constraints required to the map $g$ invoked above. It's the data of a map $s \to z$ in $\cat I$, such that $g \comp \langle p, q \rangle = \langle f \comp u, v \rangle$, which means $g \comp p = f \comp u$ and $g \comp q = v$, as desired.

	Viceversa, given $f$ and $g$ as in the previous diagram, it's easy to reconstruct $\phi$ using the universal property of the pullback.
\end{proof}

\begin{example}[Dependent lenses]
\label{ex:dlens}
	On the right hand side of~\eqref{eq:profunctor-rep} we set $\cat X = \cat Y = L$, where the latter is the functor defined in~\eqref{eq:canonical-span-action}, so that
	\begin{eqalign}
		\Optic_{L,L}&(s \nto{v} j \nfrom{v'} t, a \nto{u} i \nfrom{u'} b) :=
		\int^{i \nfrom{p} z \nto{q} j} \cat I/j(v, q_*p^* u) \times \cat I/j(q_*p^* u', v')\\
		&\iso \int^{\langle p, q \rangle : z \to i \times j} \sum_{f \in \cat I(s, a)} \cat I/(i \times j)(\langle f \comp u, v \rangle, \langle p, q \rangle) \times \cat I/j(q_*p^* u', v')\\\\
		&\iso \sum_{f \in \cat I(s, a)}  \int^{\langle p, q \rangle : z \to i \times j} \cat I/(i \times j)(\langle f \comp u, v \rangle, \langle p, q \rangle) \times \cat I/j(q_*p^* u', v')\\\\
		&\iso \sum_{f \in \cat I(s, a)} \cat I/j(v_* (f \comp u)^* u', v')
	\end{eqalign}
	One can visualize the data of this last expression as the dotted arrow in the following commutative diagram:
	\begin{diagram}
		t & {s \times_i b} && b \\
		j & s & a & i
		\arrow["{v'}"', from=1-1, to=2-1]
		\arrow["v", from=2-2, to=2-1]
		\arrow["f"', dotted, from=2-2, to=2-3]
		\arrow["u"', from=2-3, to=2-4]
		\arrow["{u'}", from=1-4, to=2-4]
		\arrow["{(f \comp u)^*u'}", from=1-2, to=2-2]
		\arrow[from=1-2, to=1-4]
		\arrow["\lrcorner"{anchor=center, pos=0.125}, draw=none, from=1-2, to=2-3]
		\arrow["{f^\sharp}"', dotted, from=1-2, to=1-1]
	\end{diagram}
	We get our usual dependent lenses when $v$ and $u$ are trivial, i.e.~when the cospans $s \nto{v} j \nfrom{v'} t$ and $a \nto{u} i \nfrom{u'} b$ are `trivial on the right'.
	To see that indeed the usual category of $\DLens$ is equivalent to this one, consider the functor sending a cospan $a \nto{u} i \nfrom{u'} b$ to the map $u^* u' : b \times_i a \to a$ obtained by pulling back $u'$ along $u$.
	On morphisms, this functor takes a diagram such as the one above and maps it to the dependent lens comprised of the middle two squares in the diagram below:
	\begin{diagram}
		t & {s \times_j t} & {s \times_a (b \times_i a)} & {b \times_i a} & b \\
		j & s & s & a & i
		\arrow["{v^*v'}"', from=1-2, to=2-2]
		\arrow[Rightarrow, no head, from=2-3, to=2-2]
		\arrow["f"', dotted, from=2-3, to=2-4]
		\arrow["u"', from=2-4, to=2-5]
		\arrow["{u'}", from=1-5, to=2-5]
		\arrow["{f^*u^*u'}", from=1-3, to=2-3]
		\arrow["\lrcorner"{anchor=center, pos=0.125}, draw=none, from=1-3, to=2-4]
		\arrow[from=1-3, to=1-4]
		\arrow[from=1-4, to=1-5]
		\arrow["{u^*u'}", from=1-4, to=2-4]
		\arrow["\lrcorner"{anchor=center, pos=0.125}, draw=none, from=1-4, to=2-5]
		\arrow["{v'}"', from=1-1, to=2-1]
		\arrow["v", from=2-2, to=2-1]
		\arrow[from=1-2, to=1-1]
		\arrow["\lrcorner"{anchor=center, pos=0.125, rotate=-90}, draw=none, from=1-2, to=2-1]
		\arrow["{f^\sharp}"', curve={height=20pt}, dotted, from=1-3, to=1-1]
		\arrow["{\exists! {f^\sharp}'}", dashed, from=1-3, to=1-2]
	\end{diagram}
	On the left, the pasting property of pullback squares allows us to conclude $f^*u^* u' = (f\comp u)^*u'$, so that the morphisms out of $s \times_i b \to s$ and $s \times_a (b \times_i a)$ are the same.
	On the right, the universal property of pullbacks tells us that the data of $f^\sharp$ such that $f^\sharp \comp v' = (f \comp u)^* u' \comp v$ is the same as that of an arrow (the dashed $\exists !$ arrow above) $s \times_a (b \times_i a) \to s \times_j t$.

	This functor is an equivalence of categories. On the one hand the action on morphisms is full and faithful, since $f$ is unaltered and $f^\sharp$ and ${f^\sharp}'$ are in bijection by universal property of the pullback.
	On the other hand it is essentially surjective since every map $p : a \to i$ is evidently isomorphic (in fact, equal) to the image of the cospan $i \equalto i \nfrom{p} a$.
\end{example}

\begin{example}[Compound optics]
	Milewski proposed compound optics \cite[§8]{milewski2022compound} as the class of optics associated to two actions of a bicategory $\twocat M$ on $\twocat M_0$-indexed categories $\cat C$ and $\cat D$ by the definition:
	\begin{equation}
	\label{eq:compound_opt}
		\Optic_{o,o'}(a,b)(s,t) := \int^{m:\twocat M(o,o')} (\Lan_{\langle m, m \rangle \action} \yo_{(a,b)})(s,t)
	\end{equation}
	Here $\yo_{(a,b)}$ denotes the Yoneda embedding of $(a,b) : \cat C(o)^\op \times \cat D(o)$; while $\langle m, m \rangle \action$ denote the action of $m:o \to o'$ in the fibrewise product of the indexed categories $\cat C^\op$ and $\cat D$, i.e.
	\begin{equation}
		\langle m, m \rangle \action : \cat C(o)^\op \times \cat D(o) \nlongto{\ \langle \cat C(m)^\op, \cat D(m) \rangle\ } \cat C(o')^\op \times \cat D(o').
	\end{equation}
	We will denote the action of $\cat C$ on morphisms of $\twocat M$ as $\action$.
	Then, using the coend formulas for Kan extensions \cite[§2]{loregian2020coend}, one has
	\begin{eqalign}
		&(\Lan_{\langle m, m \rangle \action} \yo_{(a,b)})(s,t)\\
		&\iso \int^{(x,y): \cat C(o)^\op \times \cat D(o)} \cat C(o')^\op(m \action x, s) \times \cat D(o')(m \action y, t) \times \yo_{(a,b)}(x,y)\\
		&\iso \int^{(x,y): \cat C(o)^\op \times \cat D(o)} \cat C(o')(s, m \action x) \times \cat D(o')(m \action y, t) \times \cat C(o)(x,a) \times \cat D(o)(b, y)\\
		&\iso \cat C(o')(s, m \action a) \times \cat D(o')(b \action y, t).
	\end{eqalign}
	Therefore we conclude Equation~\eqref{eq:compound_opt} and Equation~\eqref{eq:profunctor-rep} present the same object.
\end{example}

\begin{example}[Fibre optics]
\label{ex:fibre-optics}
	Fibre optics \cite{braithwaite2021fibre} are a proposed definition of a dependently-typed generalisation of optics.
	This structure is convoluted to define but makes it possible to talk about things like `stochastic dependent lenses'.
	Strictly speaking, they generalize dependent lenses by replacing the canonical action of $\Span(\cat I)$ with something satisfying the same asbtract properties.

	The data of that definition abstracts some of the properties of the actions we used to define dependent lenses above.
	Suppose we still have a locally cartesian category $\cat I$, together with an $\cat I$-indexed monoidal category of twists, $\cat M : \cat I^\op \to \MonCat$.
	Suppose, moreover, that we have indexed actegories $\cat C, \cat D : \cat I^\op \to \Actt$ such that, for every $i:\cat I$, $\cat C(i)$ and $\cat D(i)$ are $\cat M(i)$-actegories:\footnote{Here $\pi : \Actt \to \MonCat$ is the 2-fibration of actegories described in \cite[§3]{capucci2022actegories}.}
	\begin{diagram}[sep=4ex]
		& {\cat I^\op} \\
		\Actt && \MonCat
		\arrow["{\cat M}", from=1-2, to=2-3]
		\arrow["\pi"', from=2-1, to=2-3]
		\arrow["{\cat C, \cat D}"', from=1-2, to=2-1]
	\end{diagram}
	We also assume there are $\cat X_!, \cat Y_! : \cat I \to \Cat$ such that for each $f:i \to j$ in $\cat I$, $\cat X_!(f) \adj \cat C(f)$ in $\Cat$ (hence $\cat X_!$ and $\cat Y_!$ coincide with $\cat C \comp U$ and $\cat D \comp U$ on objects, where $U: \Actt \to \Cat$ forgets the actegorical structure).
	For simplicity, we denote by $f^*$ both $\cat C(f)$ and $\cat D(f)$ and by $f_!$ both $\cat X_!(f)$ and $\cat Y_!(f)$. We ask that $\cat C$ and $\cat D$ satisfy the Beck--Chevalley condition, which amounts to asking that, in the diagram below, if the square on the left is a pullback square in $\cat I$, then the square on the right is commutative in $\Cat$:
	\begin{diagram}
		\bullet & \bullet & \bullet & \bullet \\
		\bullet & \bullet & \bullet & \bullet
		\arrow["f", from=1-1, to=1-2]
		\arrow["h"', from=1-1, to=2-1]
		\arrow["k"', from=2-1, to=2-2]
		\arrow[""{name=0, anchor=center, inner sep=0}, "g", from=1-2, to=2-2]
		\arrow["\lrcorner"{anchor=center, pos=0.125}, draw=none, from=1-1, to=2-2]
		\arrow["{k_!}"', from=2-3, to=2-4]
		\arrow["{f_!}", from=1-3, to=1-4]
		\arrow[""{name=1, anchor=center, inner sep=0}, "{h^*}", from=2-3, to=1-3]
		\arrow["{g^*}"', from=2-4, to=1-4]
		\arrow[shorten <=20pt, shorten >=20pt, maps to, from=0, to=1]
	\end{diagram}
	We can use this data to define `spans with $\cat M$-apex', so that we can define the action of \emph{linear polynomials with coefficients in $\cat M$ and variables indexed by $\cat I$}.
	\begin{diagram}[sep=4ex]
		& {\cat X_z} \arrow[loop, out=60, in=120, looseness=6, "M \action_z -"'] \\
		{\cat X_i} && {\cat X_j} \\[-2ex]
		& z \arrow[loop, out=45, in=135, looseness=7, "M"', pos=0.6, squiggly]\\
		i && j
		\arrow["p"', from=3-2, to=4-1]
		\arrow["q", from=3-2, to=4-3]
		\arrow["{p^*}", from=2-1, to=1-2]
		\arrow[dotted, maps to, from=2-1, to=4-1]
		\arrow[dotted, maps to, from=2-3, to=4-3]
		\arrow["{p_!}", from=1-2, to=2-3]
		\arrow[dotted, maps to, from=1-2, to=3-2]
	\end{diagram}
	The acting bicategory is $\Mat_{\cat I}(\cat M)$, so defined:
	\begin{enumerate}
		\item objects are the same as $\cat I$,
		\item morphisms $i \to j$ are given by spans $i \nfrom{p} z \nto{q} j$ together with a choice of $M : \cat M_z$ of `apex twist',
		\item composition of $(M,p,q) : i \to j$ and $(N, r,s) : j \to k$ is given by the pink span and apex twist:
		\begin{diagram}[sep=4ex]
			&&[-6ex] \textcolor{rgb,255:red,230;green,45;blue,162}{(q^*(r)^* M) \underset{\cat M_{z \times_j w}}\otimes (r^*(q)^* N)} &[-6ex] \\
			M && \textcolor{rgb,255:red,230;green,45;blue,162}{z \times_j w} && N \\
			& \textcolor{rgb,255:red,230;green,45;blue,162}{z} && \textcolor{rgb,255:red,230;green,45;blue,162}{w} \\
			\textcolor{rgb,255:red,230;green,45;blue,162}{i} && j && \textcolor{rgb,255:red,230;green,45;blue,162}{k}
			\arrow["p"', color={rgb,255:red,230;green,45;blue,162}, from=3-2, to=4-1]
			\arrow["q", from=3-2, to=4-3]
			\arrow["r"', from=3-4, to=4-3]
			\arrow["s", color={rgb,255:red,230;green,45;blue,162}, from=3-4, to=4-5]
			\arrow["{q^*(r)}"', color={rgb,255:red,230;green,45;blue,162}, from=2-3, to=3-2]
			\arrow["{r^*(q)}", color={rgb,255:red,230;green,45;blue,162}, from=2-3, to=3-4]
			\arrow["\lrcorner"{anchor=center, pos=0.125, rotate=-45}, draw=none, from=2-3, to=4-3]
			\arrow[squiggly, from=2-1, to=3-2]
			\arrow[squiggly, from=2-5, to=3-4]
			\arrow[dotted, maps to, from=2-1, to=1-3]
			\arrow[dotted, maps to, from=2-5, to=1-3]
			\arrow[color={rgb,255:red,230;green,45;blue,162}, squiggly, from=1-3, to=2-3]
		\end{diagram}
		\item identity for an object $i : \cat I$ is given by the unit $I_i$ of $\cat M_i$,
		\item 2-morphisms $\alpha : M \twoto M' : i \to j$ are given by morphisms of $\cat M_{i \times j}$, and so are their composition and identities.
	\end{enumerate}

	Hence $\Mat_{\cat I}(\cat M)$ admits an action on the indexed category $\cat C$ (and $\cat D$) analogous to that of $\Span(\cat I)$ on the self-indexing of $\cat I$:
	\begin{diagram}
		&[-10ex]  {\Mat_{\cat I}(\cat M)} &[-10ex]&& \Cat \\[-6ex]
		& x &&& {\cat C(x)} \\
		z \arrow[loop, out=150, in=210, looseness=6, "M"', pos=0.6, squiggly] && w \arrow[loop, out=330, in=30, looseness=6, "N"', pos=0.6, squiggly] \\
		& y &&& {\cat C(y)}
		\arrow["q"', from=3-1, to=4-2]
		\arrow["p", from=3-1, to=2-2]
		\arrow["s", from=3-3, to=4-2]
		\arrow["r"', from=3-3, to=2-2]
		\arrow["\alpha", Rightarrow, from=3-1, to=3-3]
		\arrow[""{name=0, anchor=center, inner sep=0}, "q_*Mp^*"{description}, curve={height=30pt}, from=2-5, to=4-5]
		\arrow[""{name=1, anchor=center, inner sep=0}, "s_*Nr^*"{description}, curve={height=-30pt}, from=2-5, to=4-5]
		\arrow[from=1-2, to=1-5, "L_{\cat C}"]
		\arrow["{\alpha^*}", shorten <=13pt, shorten >=13pt, Rightarrow, from=0, to=1]
		\arrow[shorten <=30pt, shorten >=20pt, maps to, from=3-3, to=0]
	\end{diagram}
	Such an action is central in \cite{schreiber2014quantization}, where the action of $(M,p,q)$ is called a `pull-tensor-push operation', and it's used to abstract integral transforms.

	Then $\Optic_{L_{\cat C}, L_{\cat D}}$, in the sense of~\cref{th:tambara-rep}, is defined by:
	\begin{equation}
		\Optic_{L_{\cat C}, L_{\cat D}}((i, x,x'), (j, y,y)) := \int^{(M,p,q)} \cat C(i)(x, q_*(M \action p^*(y))) \times \cat D(j)(q_*(M \action p^*(y')), x').
	\end{equation}
	where the coend runs over $(M, i \nfrom{p} z \nto{q} j): \Mat_{\cat I}(\cat M)(i,j)$.

	Compared to definition of fibre optics in \cite{braithwaite2021fibre}, the pull and push operations are now allowed to be done along any morphism instead of just projection morphisms. That is, in the original definition the coend only ran over projection spans $i \nfrom{\pi_i} i \times j \nto{\pi_j} j$ (no restrictions on $M$ are present).
	As proven in \cite[Proposition 3.23]{schreiber2014quantization}, if we also ask for a $\cat M_!$ which provides left adjoints to each functor $\cat M(f)$, then every pull-push operation can be simulated by a pull-tensor-push along a projection span. Indeed, this is the case. in the original definition of fibre optics.
\end{example}

	\section{Conclusion}
In this note, we proposed a new definition of Tambara modules that generalizes the classical one and captures the recently proposed definitions of dependent optics \cite{vertechi2022,milewski2022compound} by a similar representation theorem as the one linking the profunctor and existential encoding of mixed optics \cite{bryce2020categoricalupdate}.

The link with double category theory is quite intriguing, and we wonder if the quick results sketched in~\cref{sec:tambara} can be put in the context of a more ample and fruitful double-categorical framework.

In particular, we feel like we ignored some double categorical information in formulating~\cref{th:tambara-rep}, which could be made not more general but more complete if we worked double-categorically throughout.
Moreover the Yoneda theory described in \cite{pare2011yoneda} should be technically useful given the observation Tambara modules are horizontal natural transformations.

Finally, the double-categorical structure should percolate to optics.
What structure does the definition of optics from doubly indexed categories impose on categories of optics? And which examples are captured by the full double categorical story which are not captured otherwise?
Indeed, we observe all the examples we have in~\cref{sec:examples} concern actions of framed bicategories \cite{shulman_framed_bicategories}, which are ambigously bi- and double categories.

	\printbibliography

@unpublished{milewski_polylens,
	author = {Bartosz Milewski},
	note = {arXiv preprint available at \url{https://bartoszmilewski.com/2021/12/07/polylens/}},
	title = {Poly{L}ens},
	year = {2021}}

@unpublished{towards_foundations_categorical_cybernetics,
	author = {Matteo Capucci and Bruno Gavranovi\'{c} and Jules Hedges and Eigil Fjeldgren Rischel},
	note = {Forthcoming in \emph{Proceedings of ACT 2021}, arXiv:2105.06332},
	title = {Towards foundations of categorical cybernetics},
	year = {2021}}

@book{myers2021book,
	author = {Myers, David Jaz},
	title = {Categorical Systems Theory},
	url = {https://github.com/DavidJaz/DynamicalSystemsBook/tree/master/book},
	urldate = {2021-05-12},
	year = {2021}}

@unpublished{myers2020double,
	author = {Myers, David Jaz},
	note = {arXiv preprint available at \url{https://arxiv.org/abs/2005.05956}},
	title = {Double Categories of Open Dynamical Systems},
	year = {2020}}

@book{loregian2020coend,
	author = {Loregian, F.},
	doi = {10.1017/9781108778657},
	edition = {first},
	month = {jul},
	note = {ISBN 9781108746120},
	publisher = {Cambridge University Press},
	series = {London Mathematical Society Lecture Note Series},
	title = {Coend Calculus},
	volume = {468},
	year = {2021}}

@article{moeller_vasilakopoulou_monoidal_grothendieck,
	author = {Joe Moeller and Christina Vasilakopoulou},
	journal = {Theory and applications of categories},
	number = {31},
	pages = {1159--1207},
	title = {Monoidal {G}rothendieck construction},
	volume = {35},
	year = {2020}}

@article{bryce2020categoricalupdate,
	author = {Clarke, Bryce and Elkins, Derek and Gibbons, Jeremy and Loregian, Fosco and Milewski, Bartosz and Pillmore, Emily and Rom\'an, Mario},
	title = {Profunctor optics: {A} categorical update},
	journal={NWPT 2019},
	pages={47},
	year = {2020}
}

@article{pastrostreet2008,
  title={Doubles for monoidal categories.},
  author={Pastro, Craig and Street, Ross},
  journal={Theory and Applications of Categories},
  volume={21},
  pages={61--75},
  year={2008},
  publisher={Mount Allison University, Department of Mathematics and Computer Science}
}

@article{tambara2006,
  title={Distributors on a tensor category},
  author={Tambara, Deisuke},
  journal={Hokkaido Mathematical Journal},
  volume={35},
  pages={379--425},
  year={2006}
}

@unpublished{spivak2019generalized,
	author = {Spivak, David I},
	note = {arXiv preprint available at \url{https://arxiv.org/1908.02202}},
	title = {Generalized Lens Categories via functors $F:\mathcal C^{\rm op} \to \mathsf{Cat}$},
	year = {2019}}

@unpublished{bakovic2009simplicial,
      title={The simplicial interpretation of bigroupoid 2-torsors},
      author={Igor Bakovic},
      year={2009},
      note={arXiv preprint availble at \url{https://arxiv.org/abs/0902.3436}}
}

@article{shulman_framed_bicategories,
	author = {Shulman, Michael},
	journal = {Theory and Applications of Categories [electronic only]},
	language = {eng},
	pages = {650-738},
	publisher = {Mount Allison University, Department of Mathematics and Computer Science, Sackville},
	title = {Framed bicategories and monoidal fibrations.},
	url = {http://eudml.org/doc/229705},
	volume = {20},
	year = {2008}}

@inproceedings{gambino2013polynomial,
	title={Polynomial functors and polynomial monads},
	author={Gambino, Nicola and Kock, Joachim},
	booktitle={Mathematical Proceedings of the Cambridge Philosophical Society},
	volume={154},
	number={1},
	pages={153--192},
	year={2013},
	organization={Cambridge University Press}}

@unpublished{capucci2022actegories,
  title={Actegories for the Working Amthematician},
  author={Capucci, Matteo and Gavranovi{\'c}, Bruno},
  note={arXiv preprint available at \url{https://arxiv.org/abs/2203.16351}},
  year={2022}
}

@unpublished{braithwaite2021fibre,
  title={Fibre optics},
  author={Braithwaite, Dylan and Capucci, Matteo and Gavranovi{\'c}, Bruno and Hedges, Jules and Rischel, Eigil Fjeldgren},
  note={arXiv preprint available at \url{https://arxiv.org/abs/2112.11145}},
  year={2021}
}

@article{moskaliuk1998double,
  title={Double categories in mathematical physics},
  author={Moskaliuk, SS and Vlassov, AT},
  journal={Ukrayins' kij Fyizichnij Zhurnal},
  volume={43},
  number={6-7},
  pages={836--841},
  year={1998}
}

@unpublished{roman2020profunctor,
  title={Profunctor optics and traversals},
  author={Rom{\'a}n, Mario},
  note={arXiv preprint available at \url{https://arxiv.org/abs/2001.08045}},
  year={2020}
}

@article{pare2011yoneda,
  title={Yoneda theory for double categories},
  author={Par{\'e}, Robert},
  journal={Theory and Applications of Categories},
  volume={25},
  number={17},
  pages={436--489},
  year={2011}
}

@unpublished{milewski2022compound,
  title={Compound Optics},
  author={Milewski, Bartosz},
  note={arXiv preprint available at \url{https://arxiv.org/abs/2203.12022}},
  year={2022}
}

@unpublished{schreiber2014quantization,
  title={Quantization via Linear homotopy types},
  author={Schreiber, Urs},
  note={arXiv preprint available at \url{https://arxiv.org/abs/1402.7041}},
  year={2014}
}

@unpublished{vertechi2022,
  author = {Vertechi, Pietro},
  note   = {arXiv preprint available at \url{https://arxiv.org/abs/2204.09547}},
  title  = {Dependent Optics},
  year   = {2022}
}

@unpublished{gavranovic2022local,
  author = {Gavranovi\'c, Bruno},
  url={https://www.brunogavranovic.com/posts/2022-02-10-optics-vs-lenses-operationally.html},
  urldate={2022-04-22},
  title  = {Optics vs Lenses, Operationally},
  year   = {2022}
}
\end{document}